\input AHTOH-E.STY

\UDC{\ttsmall
512.54.05%Алгоритмические задачи в теории групп. Проблема слов
$\scriptstyle+$512.543.53%Прямые и полные прямые произведения
$\scriptstyle+$512.543.14%Группы с конечным числом образующих
$\scriptstyle+$512.543.16%Определяющие соотношения
%$\scriptstyle+$512.543.56%Сплетения
$\scriptstyle+$512.544.7%Аппроксимация групп
$\scriptstyle+$512.552%Ассоциативные кольца и алгебры
}

\title{
Residually finite algorithmically finite groups,
their subgroups and direct products}

\author{Anton A. Klyachko\quad Ayrana K. Mongush}

\address{\myAddress\ ayranamongush@gmail.com}

\abstract{%
%\narrower
\narrower
We construct an infinite finitely generated recursively presented
residually finite algorithmically finite group~$G$ answering thereby a
question of Myasnikov and Osin.  Moreover, $G$ is ``very infinite" and
``very algorithmically finite" in the sense that $G$ contains an infinite
abelian normal subgroup while all finite Cartesian powers of $G$ are
algorithmically finite (i.e., for any positive integer $n$, there is no
algorithm which writes out an infinite sequence of pairwise different
elements of~$G^n$). We also state several related problems.
}

%%%%%%%%%%%%%%%%%%%%%%%%%%%%%%%%%%%%%%%%%%%%%%%%%%%%%%%%%%%%%%%%
\s 0.
Introduction

In [MO11], it was constructed the first example of finitely generated
recursively presented infinite group which is \emph{algorithmically
finite}, in the sense that there is no algorithm that writes out an
infinite sequence of pairwise different elements of this group. Groups
having these properties (i.e. finitely generated recursively presented
infinite and algorithmically finite) are called \emph{Dehn
monsters} in [MO11].

The Dehn monsters constructed in [MO11] also have an additional
finiteness-infiniteness property: they have infinite residually
finite homomorphic images. This led the authors of
[MO11] to the following question: {\sl Do there exist
residually finite Dehn monsters?} We answer this
question positively%
\fn{When this work have been written, we discovered that an answer
to this question is contained also in~[KhM14].} 
but ask another question.

\Question 1.
Is the direct product of two algorithmically finite groups always
algorithmically finite?

We cannot answer this question (and conjecture that the answer is
negative) but the Dehn monster constructed in this paper has
a nice property: all its finite Cartesian powers
are algorithmically finite. Even more intriguing, in our opinion, question
sounds as follows.

\Question 2.
Is the wreath product of two algorithmically finite groups always
algorithmically finite? Is it true at least that the wreath
product of a finite group (e.g., the two-element groups) and an
algorithmically finite group is always algorithmically finite?

We do not know the answer and cannot even answer the
``opposite question".

\Question 3.
Does there exist a Dehn monster $D$ such that the wreath product
$\Z_2\wr D$ is algorithmically finite?

Note, however, that the monster constructed in this paper is in some sense
similar to such a wreath product: it is a semidirect product of an
infinite elementary abelian normal subgroup and
another monster.

\proclaim{Main theorem}.
There exists an infinite finitely generated recursively presented
residually finite group $G$ such that all its finite
Cartesian powers are algorithmically finite, i.e., for any positive
integer $n$, there is no algorithm that writes out an infinite sequence
of pairwise different elements of the group~$G^n$.
\newline
Moreover, the group $G$ can be chosen containing an
infinite direct power of the cyclic groups $\Z_p$ as a normal subgroup
(where $p$ is an arbitrary given prime) and the corresponding extension
splits:  $G=H\semitimes\(\bigtimes_{i=1}^\infty\gp a_p\)$.

The question about possible subgroups of Dehn monsters deserves
a special attention. Clearly, all finitely generated subgroups of monsters
are algorithmically finite themselves; in particular, all cyclic
subgroups are finite and, therefore, e.g., all solvable finitely
generated subgroups are also finite.

\Question 4.
Which groups (or which abelian groups) can be embedded in an
algorithmically finite group?  Which groups can be embedded as
normal subgroups?

The main theorem is a corollary of the following fact about algebras,
which is of independent interest.

\Th on strongly algorithmically finite algebras.
Over any finite field, there exists an infinite finitely generated
recursively presented residually finite associative
algebra~$A$ (with unity) whose all finite Cartesian powers
are algorithmically finite, i.e., for any positive integer $n$, there is
no algorithm which writes out an infinite number of pairwise different
elements of the algebra $A^n$.
\newline
Moreover, the algebra
$A$ is generated by a finite set of nilpotent elements.%

%\fn{%
%Н  с мом деле, нетрудно добиться того, чтобы \р зр{почти все} элементы
% лгебры $A$ были нильпотентными, но мы этого не док зыв ем.
%}

Our approach is based on the ideas of
[MO11], i.e. on the application of the Golod--Shafarevich theorem.
However, our proof of the existence of Dehn monsters is
simpler than that in [MO11], in spite of the fact that we should
care about additional properties of the
monsters under construction (though, the residual finiteness
is obtained at no cost within our approach).

%%%%%%%%%%%%%%%%%%%%%%%%%%%%%%%%%%%%%%%%%%%%%%%%%%%%%%%%%%%%%%%%
\s 1. Infinite-dimensionality test

Consider the free associative algebra $F\gp X$
(with unity)
with finite basis $X$
over a field $F$. This algebra consists of polynomials
in non-commuting variables with coefficient from $F$.
We always understand the degree $\deg u$ of
a polynomial $u\in F\gp X$ as the
\emph{minimal} degree of the monomials of this polynomial.
For example, ${\deg(xy-yx+xy^{\the\year}x)=2}$.

The following convenient test for infinite-dimensionality of a
graded algebra is a corollary the well-known result of
Golod--Shafarevich [GSh64] and belongs, apparently, to M.~Ershov,
see [Er12], Corollary 2.2.

\proclaim{Infinite-dimensionality test}.
If $r_n$ is the number of elements of degree $n$ in a set $R$ consisting
of homogeneous elements of the finitely generated free associative algebra
$F\gp X$, where $r_0=r_1=0$, and the series
$$
1-|X|t+H_R(t)\:=
1-|X|t+\sum\limits_{i=2}^\infty r_it^i
%\eqno{(*)}
$$
converges to a negative number
for some $t\in(0,1)$, then the quotient algebra $A=F\gp X/(R)$ is
infinite-dimensional.

%Н м пон добятся т кже две простые леммы, перв я из которых
%(почти) очевидн ,   втор я док з н  в~[Ви65].

We need also the following obvious fact.

\Lemma 1.
If a field $F$ and a set $X$ are finite, then, for any positive integers
$n$ and $d$, there are only finitely many $n$-tuples
$
(u_{11},\dots,u_{n1}),(u_{12},\dots,u_{n2}),\dots
$
of elements of the free associative algebra $F\gp{X}$ such that, for any
different $i$ and $l$, there is $s$ such that $\deg(u_{si}-u_{sl})<d$.

\Proof
This inequality means that all tuples represent different elements of the
algebra $\(F\gp X/(X)^d\)^n$, which is obviously finite.  Here, $(X)$ is 
the ideal generated by the basis and, hence, $(X)^d$ consists of all 
polynomials of degree at least $d$.

%%%%%%%%%%%%%%%%%%%%%%%%%%%%%%%%%%%%%%%%%%%%%%%%%%%%%%%%%%%%%%%%
\s 2.
Construction of the algebra $A$

Take a positive integer $\alpha$ and a recursive everywhere defined
function $f:\N\times\N\to\N$ (see the next section for a particular choice
of $\alpha$ and $f$) and consider also a recursive enumeration
$P_1,P_2,\dots$ of all programmes with empty input and output
alphabet consisting of a finite field $F$, a
finite set $X$ containing at least two elements, and three
additional symbols: ``$+$" (plus), ``," (comma) and ``;" (semicolon).
The output sequence of each such programme is treated
as a sequence of tuples of elements of a free associative algebra
$F\gp X$: elements of each tuples are separated by commas and the tuples
are separated by semicolons;
successively written symbols from $F\sqcup X$ are treated as the product,
redundant pluses and commas are ignored. For example,
for $F=\Z_3=\{0,1,2\}$ and $X=\{x,y\}$, the sequence
$$
++xy2y+,,,,221++1+++++0xy1yyy++,,;;,;xxx2112+++yxyxy22+++
$$
is treated as four
tuples (two of them are empty, and one is incomplete):
$$
(2xy^2,2);\ ();\ ();\ (x^3+yxyxy+\dots
$$

We shall construct algebra $A$ in the form $A=F\gp X/(R)$.
The algorithm writing down the set of relators $R$ looks simple.

\medskip\noindent{\bf
Main algorithm.}
Initially, the set $R$ consists of monomials
$x^\alpha$ for all $x\in X$.  Further, on a step $k$, we run the programme
$N(k)$ (in parallel with all other running programmes) and go to the
step $k+1$.

\medskip

The programme $N(k)$ (i.e. the programme $N$ inputting a
positive integer $k$) performs the following tasks.

\medskip\noindent{\bf
Programme $N(k)$:}

\item{1.}
Run the program $P_k$ (in parallel with all other running programmes).
\item{2.}
Trace  the work of $P_k$:
when $P_k$ writes a semicolon, $N$ acts as follows:
\itemitem{a)}
interrupt (pause) the programme $P_k$;
\itemitem{b)}
check that
all tuples of elements of the algebra $F\gp X$ outputted by
$P_k$ so far have the same length $n$,
i.e. the output of $P_k$ has the form:
$$
u_{11},\dots,u_{n1};
u_{12},\dots,u_{n2};
\dots;
u_{1l},\dots,u_{nl};
\qbox{for some $u_{ij}\in F\gp X$ and some $n\in\N$};
$$
if the output does not have this form, then $N$
kills the programme $P_k$ and
terminates;
\itemitem{c)}
check whether there exists $i<l$ such that
$\deg(u_{si}-u_{sl})\ge f(n,k)$
for all $s$;
if not,
then the programme~$N$ manages the programme $P_k$ to continue
the work and
keep on tracing; if such $i<l$ is found, then $N$ proceeds with
the next step;
\itemitem{d)}
it adds all homogeneous components of the differences
$u_{si}-u_{sl}$ (for all $s\in\{1,\dots,n\}$) to the set
$R$, kills the program $P_k$ and terminates.

\bigskip

\enditem
Note that, at the every moment, there is a
finite number of programmes $P_i$ working in parallel
(at most $k$ on $k$-th step of the main algorithm) and
the same number of
copies of the programme $N$
(each copy traces one programme~$P_i$).
Moreover, each copy of the program $N$ either
\- works eternally and adds nothing to the set $R$,
\- terminates on
the step 2b) and, in this case, also adds nothing to
the set of relators $R$,
\- or terminates on the step 2d) and, in this case, adds
a finite set of homogeneous relators $w_1,w_2,\dots$ of a
large degree:  $\deg w_i\ge f(n,k)$
(where $k$ is the number of this
copy of $N$); the number of added relators of the each particular degree
is at most $n$.  More precisely, we have the inequality:
$$
r_i(k)\le\cases{
0, &    if $i<f(n(k),k)$;\cr
n(k), & if $i\ge f(n(k),k)$;
}
$$
where $r_i(k)$ is the number of relators of degree $i$ added
by the $k$-th copy of the programme $N$ and $n(k)$ is a
the length of tuples
outputted by $P_k$ (if $P_k$ writes out tuples of different
lengths, or some incorrect output, or does not write anything, then we
assume $n(k)=\infty$).

%%%%%%%%%%%%%%%%%%%%%%%%%%%%%%%%%%%%%%%%%%%%%%%%%%%%%%%%%%%%%%%%
\s 3.
Infinite-dimensionality of the algebra $A$

To apply the infinite-dimensionality test, we have to
estimate the sum of the Golod--Shafarevich series:
$$
\eqalign{
H_R(t)\:=&
\sum\limits_{i=2}^\infty r_it^i
=
|X|t^\alpha+\sum_{k=1}^\infty\(\sum\limits_{i=2}^\infty r_i(k)t^i\)
\le
|X|t^\alpha+\sum_{k=1}^\infty\(\sum\limits_{i=f(n(k),k)}^\infty n(k)t^i\)
=
\cr
=&
|X|t^\alpha+\sum_{k=1}^\infty\(t^{f(n(k),k)}n(k){1\over1-t}\)
=
|X|t^\alpha+{1\over1-t}\sum_{k=1}^\infty\(t^{f(n(k),k)}n(k)\).
}
$$
Setting $t={1\over 2}$ and taking into account that $x<2^x$ for $x\in\N$,
we obtain
$$
H_R\({1\over2}\)=
{|X|\over 2^\alpha}+2\sum_{k=1}^\infty\({n(k)\over2^{f(n(k),k)}}\)
<
{1\over 2^{\alpha-|X|}}+2\sum_{k=1}^\infty\({1\over2^{f(n(k),k)-n(k)}}\).
$$
Now, put $f(n,k)=n+k+2$ and $\alpha=|X|+1$ and obtain
$$
H_R\({1\over2}\)
<
{1\over 2}+2\sum_{k=1}^\infty\({1\over4\cdot2^k}\)=1,
\qqbox{i.e.}
1-{1\over2}|X|+H_R\({1\over2}\)<0
\qbox{if $|X|\ge4$}.
$$
According to the test from Section 1, the graded algebra
$A=F\gp X/(R)$ is infinite-dimensional. Certainly, this algebra is
residually finite because any finitely generated graded algebra over
finite field is approximated by its finite quotient algebras $A/(X)^n$.

%%%%%%%%%%%%%%%%%%%%%%%%%%%%%%%%%%%%%%%%%%%%%%%%%%%%%%%%%%%%%%%%
\s 4.
Algorithmic finiteness of Cartesian powers of $A$

Suppose that there is a programme $P$ writing out an infinite sequence of
pairwise different elements of the algebra~$A^n$. Certainly, this
programme can be
assumed to have the output alphabet
$X\sqcup F\sqcup\{``{+}",``,",``;"\}$
and writes out pairwise different elements
of the algebra $A^n$ in the prescribed format:
$$
u_{11},\dots,u_{n1};
u_{12},\dots,u_{n2};
\dots,
\qbox{where $u_{ij}\in F(X)$}
$$
(any programme can be transformed into this form).

The programme $P$ is assigned a number
$k$ by our enumeration of programmes, i.e. $P=P_k$. Then
the programe~$N(k)$ is launched on $k$-th step of our main algorithm
(in parallel with other working
programmes); $N(k)$, in its turn, launches the programme $P_k=P$
and trace it. Two cases are possible.

{\noindent\bf
Case I:}
$\deg(u_{si}-u_{sl})\ge f(n,k)$
for some distinct $i$ and $l$ and all $s$.
We assume that $i<l$ and $l$ is the minimal number with these properties.
Then, after writing out the $l$-th semicolon, the programme $P=P_k$ is
interrupted on step 2a) of the tracer-programme $N(k)$. Further, on steps
2b) and 2c), the verification is completed successfully and, on step~2d),
all homogeneous components of the differences ${u_{si}-u_{sl}}$ are added
to the set $R$ (for all $s\in\{1,\dots,n\}$). This means that the tuples
$(u_{1i},\dots,u_{ni})$ and $(u_{1l},\dots,u_{nl})$ outputted by the
programme~$P$ represent the same element of the algebra $A^n$ and we come
to a contradiction.

{\noindent\bf
Case II:}
for any distinct $i$ and $l$, there is $s$ such that
$\deg(u_{si}-u_{sl})< f(n,k)$. This contradicts Lemma 1.

\medskip
The theorem on strongly algorithmically finite algebras is proven.

%%%%%%%%%%%%%%%%%%%%%%%%%%%%%%%%%%%%%%%%%%%%%%%%%%%%%%%%%%%%%%%%
\s 5.
Proof of the main theorem

Take the algebra $A$ over a residue field $\Z_p$ with finite
generating set $X$ from the theorem on
strongly algorithmically finite algebras and consider the set of matrices
$$
G=\pmatrix{
H&A\cr
0&1
},
%\left\{
%\right\}
$$
where $H$ is the subgroup of the multiplicative group of $A$ generated
by all
elements of the form $1+x$, where $x\in X$ (these elements are invertible,
because elements of $X$ are nilpotent). Clearly, $G$ is a group:
$$
\pmatrix{
h&a\cr
0&1
}
\pmatrix{
h'&a'\cr
0&1
}
=
\pmatrix{
hh'&a+ha'\cr
0&1
},
$$
and $G$ is a semidirect product of
$$
\hbox{the normal subgroup }
\pmatrix{
1&A\cr
0&1
}\iso\bigoplus\limits_{i=1}^\infty\Z_p
\qbox{
and the non-normal subgroup }
\pmatrix{
H&0\cr
0&1
}\iso H.
$$
The group $G$ is finitely generated (($|X|+1$)-generated) as
it is generated by  matrices
$
\pmatrix{
1+X&0\cr
0&1
}
$
and
$
\pmatrix{
1&1\cr
0&1
}
$,
because, obviously, elements of the set $1+X$ generate
$A$ as a ring.
All Cartesian powers $G^n$ are algorithmically finite because
all Cartesian powers of the algebra $A$ are algorithmically finite.
Clearly, the group $G$ is recursively presented since the
algebra $A$ is recursively presented.
The residual finiteness of the group $G$ also follows immediately from
the residual finiteness of the algebra $A$,
because the group of invertible matrices (and the algebra of all matrices)
over residually finite algebra is residually finite.
This completes the proof of the main theorem.

%%%%%%%%%%%%%%%%%%%%%%%%%%%%%%%%%%%%%%%%%%%%%%%%%%%%%%%%%%%%%%%%
\REFERENCES

%\[Ви65]
%Э. Б. Винберг.
%К теореме о бесконечномерности  ссоци тивной  лгебры //
%Изв. АН СССР. Сер. м тем., 29:1 (1965), 209--214.

\[GSh64]
E. S. Golod, I. R. Shafarevich.
On the class field tower//
Izv. Acad. Nauk SSSR. Ser. Mat., 28:2 (1964), 261--272.

\[Er12]
M. Ershov.
Golod--Shafarevich groups: a survey //
Int. J. Algebra Comput. 22:5 (2012), 1230001, 68 pp.
\newline
See also arXiv:1206.0490.

\[KhM14]
B. Khoussainov, A. Miasnikov
Finitely presented expansions of groups, semigroups, and algebras //
Trans. Amer. Math. Soc. 366 (2014), 1455-1474

\[MO11]
A. Myasnikov, D. Osin.
Algorithmically finite groups //
J. Pure Appl. Algebra 215:11 (2011), 2789--2796.
\newline
See also arXiv:1012.1653.

\end